%%% GILKEY/IVANOVA OSSERMAN CONJECTURE
\input amstex
\documentstyle{amsppt}

\hsize=124mm \vsize=186mm
\NoRunningHeads%\NoBlackBoxes
%\newcount\qcts\newcount\qcta\qcta=1\newcount\qcteq\newcount\qcthead
%\def\setref#1{\global\qcteq=0\global\advance\qcts by\number\qcta
%       \immediate\write3{
%    \string\def\string#1\string{\number\qcthead.\number\qcts\string}}}
%\def\sethead#1{\global\qcts=0\global\advance\qcthead by\number\qcta
%       \immediate\write3{
%    \string\def\string#1\string{\number\qcthead\string}}}
%\def\seteqn#1{\global\advance\qcteq by\number\qcta\immediate\write3{
%    \string\def\string#1\string{\number\qcthead.\number\qcts.\eqnpbg\string}}}
%\def\eqnpbg{\ifnum\qcteq=1 a\else\ifnum\qcteq=2 b\else
%   \ifnum\qcteq=3 c\else\ifnum\qcteq=4 d\else
%   \ifnum\qcteq=5 e\else\ifnum\qcteq=6 f\else\ifnum\qcteq=7 g\else
%   \ifnum\qcteq=8 h\else\ifnum\qcteq=9 i\else
%   \ifnum\qcteq=10 j\else\ifnum\qcteq=11 k\else\ifnum\qcteq=12 l\else
%   \ifnum\qcteq=13 m\else\ifnum\qcteq=14 n\else\ifnum\qcteq=15 o\else
%   \ifnum\qcteq=16 p\else\ifnum\qcteq=17 q\else
%   \ifnum\qcteq=18 r\else\ifnum\qcteq=19 s\else\ifnum\qcteq=20 t\else
%   \ifnum\qcteq=21u\else\ifnum\qcteq=22 v\else
%   \ifnum\qcteq=23 w\else\ifnum\qcteq=24 x\else
%   \ifnum\qcteq=25 y\else\ifnum\qcteq=26 z\else *
%    \fi\fi\fi\fi\fi\fi\fi\fi\fi\fi\fi\fi\fi\fi\fi\fi\fi\fi\fi\fi\fi\fi\fi\fi\fi\fi}
%\newcount\qct\newcount\qcta\qct=0\qcta=1
%\def\pbgkey#1{\key{#1}\global\advance\qct by\number\qcta
%       \immediate\write3{\string\def\string#1\string{\number\qct\string}}}
\def\pcases#1#2{\hbox{#1}&\hbox{if\
   #2}}

%\immediate\input W04macro.tex            %%% Bibliographic references & Headings
%\immediate\openout 3=W04macro.tex        %%% Bibliographic references file
\def\nmonth{\ifcase\month\ \or January\or
   February\or March\or April\or May\or June\or July\or August\or
   September\or October\or November\else December\fi}
\def\range{\operatorname{range}}
\def\Pspan{\operatorname{span}}
\def\rank{\operatorname{rank}}
\def\id{\operatorname{Id}}

\font\pbglie=eufm10
\def\RR{\text{\pbglie R}}
\def\real{{\Bbb{R}}}

\def\qctSA{1}
 \def\arefa{1.1}
 \def\arefaa{1.1.a}
 \def\arefac{1.1.b}
 \def\qctSB{2}
 \def\brefa{2.1}
 \def\brefb{2.2}
 \def\brefba{2.2.a}
 \def\brefbc{2.3}
 \def\CREFbba{2.3.a}
 \def\brefd{2.4}
 \def\brefe{2.5}
 \def\breff{2.6}
 \def\brefg{2.7}
 \def\qctC{3}
 \def\crefa{3.1}
 \def\crefaa{3.1.a}
 \def\qctD{4}
 \def\drefa{4.1}
 \def\drefb{4.2}
 \def\drefc{4.3}
 \def\drefe{4.4}
\def\refBBG{1}
\def\refBBGZ{2}
\def\refBCGHV{3}
\def\refChi{4}
\def\refDD{5}
\def\refGKV{6}
\def\refGVV{7}
\def\refG{8}
\def\refGa{9}
\def\refGb{10}
\def\refGiIv{11}
\def\refGiStavrov{12}
\def\refGSVx{13}
\def\refGSV{14}
\def\refO{15}
\def\refSa{16}
\def\refSt{17}
\def\refSV{18}
\def\sethead#1{}
\def\seteqn#1{}
\def\setref#1{}
\def\pbgkey#1{\key{#1}}
\topmatter
\title The Jordan normal form of higher order Osserman algebraic curvature
tensors\endtitle 
\author Peter Gilkey\dag\ and Raina Ivanova\ddag\endauthor
\address PG: Mat.  Dept., University of Oregon, Eugene Oregon
97403 USA\newline\phantom{http:}
   http://darkwing.uoregon.edu/$\sim$gilkey\endaddress
\email gilkey\@darkwing.uoregon.edu\endemail
\address RI: Dept. of Descriptive Geometry, University of Architecture,
Civil\newline\phantom{http:} Engineering \& Geodesy, 1, Christo Smirnenski Blvd., 1421 Sofia,
Bulgaria\newline\phantom{http:}
   http://hopf.uoregon.edu/$\sim$ivanovar\endaddress
\email ivanovar\@is.tsukuba.ac.jp and ivanovar\@hopf.uoregon.edu\endemail
\thanks 2000 {\it Mathematics Subject Classification.} 53B20\endthanks
\thanks {\it Key words and phrases.} 
Higher order Jacobi operator, Osserman algebraic curvature tensors, Jordan
Osserman algebraic curvature tensors\endthanks
\thanks\dag Research partially supported by the NSF (USA) and the MPI (Leipzig)
\endthanks
\thanks\ddag Research partially supported by JSPS Post Doctoral Fellowship Program (Japan)
and the MPI (Leipzig)
\endthanks
%\subjclass Primary 53B20\endsubjclass
%\keywords Higher order Jacobi operator, Osserman algebraic curvature tensors, Jordan
%Osserman algebraic curvature tensors\endkeywords
\abstract We construct new examples of algebraic curvature tensors so that the Jordan
normal form of the higher order Jacobi operator is constant on the Grassmannian of subspaces
of type $(r,s)$ in a vector space of signature $(p,q)$. We then use these examples to establish
some results concerning higher order Osserman and higher order Jordan Osserman algebraic
curvature tensors. \endabstract
\endtopmatter
\def\mathcal#1{{\Cal{#1}}}

\def\JJ{{\mathcal{J}}}

\font\pbglie=eufm10

\def\Gr{\text{Gr}}

\sethead\qctSA
\head\S\qctSA\ Introduction\endhead

A $4$ tensor $R$ is said to be an algebraic curvature
tensor if it satisfies the well known symmetries of the Riemannian curvature tensor, i.e.
$$\eqalign{
&R(x,y,z,w)=-R(y,x,z,w)=R(z,w,x,y),\quad\text{ and}\cr
&R(x,y,z,w)+R(y,z,x,w)+R(z,x,y,w)=0.
}$$
It is clear that the Riemann curvature tensor defines an algebraic curvature tensor at each
point of the manifold. Conversely, every algebraic curvature tensor is geometrically
realizable \cite{\refGb}. We remark that it is often convenient to study certain geometric
problems in a purely algebraic setting.

Let $R$ be an algebraic curvature tensor on a vector space $V$ of signature $(p,q)$. The {\it
Jacobi operator} $\JJ_R$ is the self-adjoint linear map defined by:
$$\JJ_R(v)y:=R(y,v)v.$$
Here the natural domains of definition are the {\it pseudo-spheres}
of unit timelike ($-$) and spacelike $(+)$ vectors in $V$:
$$S^\pm(v):=\{v\in V:(v,v)=\pm1\}.$$

Let $(M,g)$ be a pseudo-Riemannian manifold of signature $(p,q)$ and let ${}^gR$ be the
curvature tensor of the Levi-Civita connection. If $(M,g)$ is Riemannian (i.e. $p=0$), and if it
is flat or it is a local rank
$1$ symmetric space, then the set of local isometries acts transitively on the unit sphere
bundle
$S(M,g)$. Consequently, the eigenvalues of
$\JJ_{{}^gR}$ are constant on $S(M,g)$.  Osserman
\cite{\refO} wondered if the converse holds; later authors called this problem the {\it
Osserman conjecture}. The conjecture has been established by Chi \cite{\refChi} for
Riemannian manifolds of dimension
$m$, where $m=4$, where $m$ is odd, or where $m\equiv2$ mod $4$. However, it is known
\cite{\refG} that there exist Riemannian Osserman algebraic curvature tensors which are not flat
and which are not the curvature tensors of rank $1$ symmetric spaces. We also refer to
\cite{\refDD,\refGSV} for related results.   

In any signature, we say that an algebraic curvature tensor is {\it Osserman} if the
eigenvalues of
$\JJ_R$ are constant on $S^\pm(V)$. Similarly, we say that a
pseudo-Riemannian manifold
$(M,g)$ is {\it Osserman} if the eigenvalues of $\JJ_{{}^gR}$ are constant on the pseudo-sphere
bundles  $S^\pm(M,g)$. 

In the Lorentzian setting ($p=1$), it is known
\cite{\refBBG,\refGKV} that an Osserman algebraic
curvature tensor has constant sectional curvature. Thus we may draw the geometric consequence
that a Lorentzian Osserman manifold has constant sectional curvature; the geometry of
such manifolds is very special. 

In higher signatures, although there are some partial
results known, the classification is far from complete. In
particular, it is known that there exist pseudo-Riemannian Osserman manifolds which are neither
flat nor local rank
$1$ symmetric spaces \cite{\refBBGZ,\refBCGHV,\refGVV}. 

In the Riemannian setting, any self-adjoint linear map is diagonalizable; thus the eigenvalues
determine the Jordan normal form (i.e. the conjugacy class). However, this is not true in
higher signature so we have to differentiate between the eigenvalue structure and the Jordan
normal form. We say that an algebraic curvature tensor is {\it Jordan Osserman} if
the Jordan normal form of
$\JJ_R$ is constant on
$S^\pm(V)$. There exist algebraic curvature tensors which are Osserman but which are not
Jordan Osserman. Furthermore, there exist Jordan Osserman algebraic curvature tensors whose
Jordan normal form is arbitrarily complicated
\cite{\refGiIv}.

The Jacobi operator was originally defined on $S^\pm(V)$. However, since we have
$\JJ_R(tv)=t^2\JJ_R(v)$, we can also regard
$\JJ_R$ as being defined on the projective spaces of non-degenerate lines in $V$ by
setting
\setref\arefa
\seteqn\arefaa
$$\JJ_R(\text{span}\{v\}):=(v,v)^{-1}\JJ_R(v)\quad\text{ if }\quad(v,v)\ne0.
\tag\arefaa$$
 Stanilov has extended $\JJ_R$ to
non-degenerate subspaces of arbitrary dimension. Let $\sigma$ be a
non-degenerate subspace of $V$. If
$\{v_i\}$ is a basis for $\sigma$, then let $h_{ij}:=(v_i,v_j)$ describe the restriction of the
metric on
$V$ to the subspace $\sigma$. Since $\sigma$ is non-degenerate, the matrix $(h_{ij})$ is
invertible and we let
$(h^{ij})$ be the inverse. The {\it higher order Jacobi operator} is
defined
\cite{\refSa} by generalizing equation (\arefaa):
$$\JJ_R(\sigma)y:=\textstyle\sum_{ij}h^{ij}R(y,v_i)v_j;$$
it is independent of the basis chosen.
This extends the natural domains of $\JJ_R$ to the {\it Grassmannians} $\Gr_{r,s}(V)$ of
subspaces of
$V$ which have signature $(r,s)$. Let
$\{e_1,...,e_{r+s}\}$ be an orthonormal basis for
$\sigma\in\Gr_{r,s}(V)$. Let $\varepsilon_i:=(v_i,v_i)$. Then we can express $\JJ_R$ more
simply as:
\seteqn\arefac
$$\JJ_R(\sigma)=\varepsilon_1\JJ_R(e_1)+...+\varepsilon_k\JJ_R(e_k).\tag\arefac$$

Now we extend the notions `Osserman' and `Jordan Osserman' to the higher order context. We say
that $R$ is {\it Osserman of type $(r,s)$} if the eigenvalues of $\JJ_R(\cdot)$ are constant on
$\Gr_{r,s}(V)$. Furthermore, $R$ is said to be {\it Jordan Osserman of type $(r,s)$} if the
Jordan normal form of $\JJ_R(\cdot)$ is constant on $\Gr_{r,s}(V)$. Since the Jordan normal
form determines the eigenvalues, it is immediate that if
$R$ is Jordan Osserman of type $(r,s)$, then $R$ is Osserman of type $(r,s)$; the reverse
implication can fail - see Remark \brefg.

We say that a pair $(r,s)$ is {\it admissible} if $\Gr_{r,s}(V)$ is non-empty and does not
consist of a single point. Equivalently, this means that:
$$0\le r\le p,\quad 0\le s\le q,\quad\text{ and }\quad1\le r+s\le\dim V-1.$$

In Section \qctSB, we state the main Theorems of this paper concerning the higher order
Jacobi operator, Osserman algebraic curvature tensors, and Jordan Osserman algebraic curvature
tensors. In Theorem
\brefa, we summarize previously known results for Osserman algebraic curvature
tensors. Theorem
\brefb\ deals with Jordan Osserman duality. In Theorem \brefbc, we present examples
due to \cite{\refGSVx} and note that previously known results for these examples can be extended
from the Osserman to the Jordan Osserman setting. In Theorem
\brefd, we construct new examples of algebraic curvature tensors which are Jordan
Osserman for certain but not all values of
$(r,s)$. We use these examples to draw certain conclusions about the relationship between the
various concepts which we have introduced.
In Section \qctC, we prove Theorem
\brefa. In Section \qctD, we prove Theorem \brefd.

\sethead\qctSB
\head\S\qctSB\ Statement of results\endhead

In the following theorem, we summarize previously known facts concerning Osserman algebraic
curvature tensors (see \cite{\refGSVx} for details).  Assertion
(2) is a duality result. Moreover, assertion (3) shows that when we consider the
eigenvalue structure, only the value of $r+s$ is relevant. Consequently, we shall say that
$R$ is {\it $k$-Osserman} if $R$ is Osserman of type $(r,s)$ for any, and hence for all,
admissible pairs
$(r,s)$ for which
$k=r+s$. In particular, if $p>0$ and if $q>0$, then the eigenvalues of $\JJ_R$ are constant on
$S^+(V)$ if and only if they are constant on $S^-(V)$. 
\setref\brefa
\proclaim{\brefa\ Theorem} Let $R$ be an algebraic curvature tensor on a vector space $V$ of
signature $(p,q)$. Let $R$ be Osserman of type $(r,s)$, where $(r,s)$ is an admissible pair. Then
we have:\roster
\smallskip\item $R$ is Einstein.
\smallskip\item $R$ is Osserman of type $(p-r,q-s)$. 
\smallskip\item If $(\tilde r,\tilde s)$ is an admissible pair with $r+s=\tilde r+\tilde s$, then
$R$ is Osserman of type $(\tilde r,\tilde s)$.\endroster\endproclaim

Instead of the eigenvalue structure, we can consider the Jordan normal form and establish
a similar duality result.

\setref\brefb
\proclaim{\brefb\ Theorem} Let $R$ be an algebraic curvature tensor on a vector space $V$ of
signature
$(p,q)$. Let $(r,s)$ be an admissible pair. If $R$ is Jordan Osserman of type $(r,s)$, then $R$
is Jordan Osserman of type
$(p-r,q-s)$.\endproclaim

The only examples given in the literature \cite{\refGSVx} may be described as follows. Let
$$R_{\id}(x,y)z:=(y,z)x-(x,z)y$$
denote the algebraic curvature tensor of constant sectional curvature. If $\phi$ is a
skew-adjoint map of $V$, then we may define:
$$R_\phi(x,y)z:=(\phi y,z)\phi x-(\phi x,z)\phi y-2(\phi x,y)\phi z.$$
We showed \cite{\refGiIv} that $R_\phi$ is an algebraic curvature tensor. We then have
\seteqn\brefba
$$\eqalign{
&\JJ_{R_{\id}}(x)y=(x,x)y-(x,y)x\quad\text{and}\cr
&\JJ_{R_\phi}(x)y=3(\phi x,y)\phi x.}\tag\brefba$$
Let $\phi$ be a skew-adjoint map with $\phi^2=\pm\id$. Let $c_0$ and $c_1$ be real
constants. We set
$R:=c_0R_{\id}+c_1R_\phi$. If
$\sigma$ is a non-degenerate subspace, then $\JJ_R(\sigma)$ is diagonalizable; thus the
eigenvalue structure determines the Jordan normal form. The following result in the Jordan
Osserman context then follows from the corresponding result in the Osserman context
\cite{\refGSVx}.

\setref\brefbc
\proclaim{\brefbc\ Theorem} Let $V$ be a vector space of signature $(p,q)$. Let $\phi$ be a
skew-adjoint map of $V$ with $\phi^2=\pm\id$.\roster
\smallskip\item The algebraic curvature tensors $R_{\id}$ and $R_\phi$
are Jordan Osserman of type $(r,s)$ for every admissible pair $(r,s)$.
\smallskip\item Let $c_0$ and $c_1$ be non-zero constants. Let $R=c_0R_{\id}+c_1R_\phi$. Then
$R$ is Jordan Osserman of types $(1,0)$, $(0,1)$, $(p-1,q)$, and $(p,q-1)$. Furthermore,
$R$ is not Osserman of type $(r,s)$ for other values of $(r,s)$.
\endroster
\endproclaim

As noted above, the operators $\JJ_R(\sigma)$ associated to the algebraic curvature tensors
discussed in Theorem
\brefbc\ are all diagonalizable. We now construct algebraic curvature tensors so
$\JJ_R(\sigma)$ has non-trivial Jordan normal form. Let $p\ge2$ and let
$q\ge2$. Let
$\{e_1^-,...,e_p^-,e_1^+,...,e_q^+\}$ be an orthonormal basis for $V$, where the vectors
$\{e_1^-,...,e_p^-\}$ are timelike and the vectors $\{e_1^+,...,e_q^+\}$ are spacelike. 
Let $a$ be a positive integer with $2a\le\min(p,q)$. We define a skew-adjoint linear map
$\Phi_a$ of
$V$ by setting:
\seteqn\CREFbba
$$\phi_ae_k^\pm=\cases
\pcases{$\pm(e_{2i}^-+e_{2i}^+)$}{$k=2i-1\le2a$,}\cr
\pcases{$\mp(e_{2i-1}^-+e_{2i-1}^+)$}{$k=2i\le2a$,}\cr
\pcases{$0$}{$k>2a$.}\endcases\tag\CREFbba$$
The map $\Phi_a$ is the direct sum of $a$ different $4\times 4$
`blocks'; the subspaces  spanned by $\{e_{2i-1}^-,e_{2i}^-,e_{2i-1}^+,e_{2i}^+\}$ are
invariant under the action of $\Phi_a$ for $1\le i\le a$.

\medbreak We can interchange the roles of spacelike and timelike vectors by changing the sign of
the inner product. Thus we may always assume that $p\le q$. The following result giving new
families of examples is in many ways the main result of this paper.

\setref\brefd
\proclaim{\brefd\ Theorem} Let $\Phi_a$ be the skew-adjoint linear  map on the vector
space $V$ of signature $(p,q)$ which is defined in equation
{\rm(\CREFbba)}. Let $R_a$ be the associated curvature tensor. Assume $p\le q$. We
have:
\roster
\smallskip\item $R_a$ is $k$ Osserman for $1\le k\le\dim V-1$.
\smallskip\item Suppose that $2a<p$. Then $R_a$ is Jordan Osserman of type $(p,0)$ and $(0,q)$;
$R_a$ is not Jordan Osserman of type $(r,s)$ otherwise.
\smallskip\item Suppose that $2a=p<q$. Then $R_a$ is Jordan Osserman of type $(r,0)$ and of type
$(r,q)$ for any $1\le r\le p-1$; $R_a$ is not Jordan Osserman otherwise.
\smallskip\item Suppose that $2a=p=q$. Then $R_a$ is Jordan Osserman of type $(r,0)$, of type
$(r,q)$, of type $(0,s)$, and of type $(p,s)$ for $1\le r\le p-1$ and $1\le s\le q-1$; $R_a$ is not Jordan
Osserman otherwise.
\endroster
\endproclaim

\setref\brefe
\definition{\brefe\ Remark} We may use assertion (1) of Theorem \brefd\
to see that assertion (2) of Theorem
\brefa\ does not generalize to the Jordan Osserman context; there exist algebraic
curvature tensors which are $k$ Osserman for all $k$, which are Jordan Osserman of type $(p,0)$
and
$(0,q)$, and which are not Jordan Osserman of type $(r,s)$ for other values of $(r,s)$. Thus we
can not determine whether or not $R$ is Jordan Osserman of type $(r,s)$ only from
$k=r+s$.\enddefinition

\setref\breff
\definition{\breff\ Remark}
We suppose $2\le k\le\dim V-2$ to ensure that we are truely in the higher order setting.
The
$k$ Osserman algebraic curvature tensors have been classified in the Riemannian and in the
Lorentzian settings
\cite{\refGa,\refGiStavrov}; all these curvature tensors have constant sectional curvature and
hence are $1$ Osserman. Again, we may use assertion (1) of Theorem
\brefd\ to see that a similar assertion fails for a Jordan Osserman algebraic curvature
tensor in the higher signature setting.\enddefinition

\setref\brefg
\definition{\brefg\ Remark} The algebraic curvature tensors described in Theorem
\brefd\ show that Osserman of type $(r,s)$ does not imply Jordan Osserman of type $(r,s)$.
\enddefinition

\medbreak It is useful to give a graphical representation of Theorem \brefd. We may
think of the values of $(r,s)$ as the points with integer coordinates in the rectangle 
$$\RR:=\{(r,s):0\le r\le p,\ 0\le s\le q\}.$$
The two corners $(0,0)$ and $(p,q)$ of $\RR$ are
excluded as inadmissible; $R_a$ is always Jordan Osserman at the other two corners $(p,0)$ and
$(0,q)$. These two corner points are the only values for which $R_a$ is Jordan Osserman if
$2a<p\le q$. If $2a=p<q$, then $R_a$ is Jordan Osserman on the two edges of $\RR$ parallel to
the $r$ axis. If $2a=p=q$, then $R_a$ is Jordan Osserman on the boundary of $\RR$. The values
for which
$R_a$ is Jordan Osserman are graphically represented by the three different pictures given
below. Entries with `$\star$' are points where $R$ is Jordan Osserman, entries with `$\circ$'
are points where $R$ is not Jordan Osserman, and entries with `$-$' are inadmissible points.
The $r$-axis is horizontal and the $s$-axis is vertical.
\medbreak\noindent\centerline{$2a<p\le q$
\hglue 1.8cm$2a=p<q$\hglue 1.8cm$2a=p=q$}
\smallbreak\noindent\centerline{\hbox{\vrule\vbox{\offinterlineskip
\hrule \halign {\phantom{.}#\quad&#\quad&#\quad&#\quad&#\phantom{.}\phantom{\vrule height
12pt}\cr
$\star$&$\circ$&$...$&$\circ$&$-$\cr
$\circ$&$\circ$&$...$&$\circ$&$\circ$\cr
$...$&$...$&$...$&$...$&$...$\cr
$\circ$&$\circ$&$...$&$\circ$&$\circ$\cr
$-$    &$\circ$&$...$&$\circ$&$\star$
\cr\phantom{\vrule
height 4pt}\cr}\hrule}\vrule}
\hbox{\vrule\vbox{\offinterlineskip
\hrule \halign {\phantom{.}$#$\quad&$#$\quad&$#$\quad&$#$\quad&$#$\phantom{.}\phantom{\vrule
height 12pt}\cr
\star&\star&...&\star&-\cr
\circ&\circ&...&\circ&\circ\cr
...&...&...&...&...\cr
\circ&\circ&...&\circ&\circ\cr
-    &\star&...&\star&\star
\cr\phantom{\vrule
height 4pt}\cr}\hrule}\vrule\ 
\hbox{\vrule\vbox{\offinterlineskip
\hrule \halign {\phantom{.}$#$\quad&$#$\quad&$#$\quad&$#$\quad&$#$\phantom{.}\phantom{\vrule
height 12pt}\cr
\star&\star&...&\star&-\cr
\star&\circ&...&\circ&\star\cr
...&...&...&...&...\cr
\star&\circ&...&\circ&\star\cr
-    &\star&...&\star&\star
\cr\phantom{\vrule
height 4pt}\cr}\hrule}\vrule}}}

\sethead\qctC
\head\S\qctC\ Jordan Osserman duality\endhead

Theorem \brefa\ (2) is a duality result for Osserman algebraic curvature tensors which
was proved in \cite{\refGSVx}. We generalize that proof to establish the corresponding
duality result for Jordan Osserman algebraic curvature tensors given in Theorem
\brefb. 

Let $R$ be Jordan Osserman of type $(r,s)$ on a vector space $V$ of signature $(p,q)$. We must
show that $R$ is Jordan Osserman of type
$(p-r,q-s)$, i.e. that the Jordan normal form of $\JJ_R(\cdot)$ is constant on
$\Gr_{p-r,q-s}(V)$. 

By Theorem \brefa\ (1), $R$ is Einstein. Let $\{e_1,...,e_{p+q}\}$ be any
orthonormal basis for $V$. If $\varepsilon_i:=(e_i,e_i)$, then
$\textstyle\sum_i\varepsilon_iR(y,e_i,e_i,x)=c(y,x)$, where $c$ is the Einstein constant.
This implies that:
\setref\crefa
\seteqn\crefaa
$$\textstyle\sum_i\varepsilon_i\JJ_R(e_i)=c\cdot\id.\tag\crefaa$$

Let $\tau\in\Gr_{p-r,p-s}(V)$. Let $\sigma:=\tau^\perp\in\Gr_{r,s}(V)$ be the orthogonal
complement of $\tau$. We construct an adapted orthonormal basis for $V$ as follows. Let
$\{e_1,...,e_{r+s}\}$ be an orthonormal basis for
$\sigma$ and let
$\{e_{r+s+1},...,e_{p+q}\}$ be an orthonormal basis for $\tau$. We use equations (\arefac)
and (\crefaa) to see that
$$\JJ_R(\sigma)+\JJ_R(\tau)=\textstyle\sum_{1\le
i\le r+s}\varepsilon_i\JJ_R(e_i)+\textstyle\sum_{r+s+1\le i\le p+q}\varepsilon_i
\JJ_R(e_i)=c\cdot\id.$$
Thus the Jordan normal form of $\JJ_R(\tau)$ is determined by the Jordan normal form of
$\JJ_R(\sigma)$. By hypothesis the Jordan normal form of $\JJ_R(\sigma)$ is constant on
$\Gr_{r,s}(V)$. Thus we may conclude that the Jordan normal form of $\JJ_R(\tau)$ is
also constant on $\Gr_{p-r,p-s}(V)$. \qed

\sethead\qctD
\head\S\qctD\ Examples of Jordan Osserman algebraic curvature tensors\endhead

Throughout this section, we shall let $V$ be a vector space of signature $(p,q)$, we shall let
$\Phi_a$ be the skew-adjoint linear map on $V$ which is defined in equation
{\rm(\CREFbba)}, we shall let $R_a$ be the associated curvature
tensor, and we shall let $\JJ_a$ be the associated Jacobi operator.
We begin the proof of Theorem \brefd\ with the following
observation:

\setref\drefa
\proclaim{\drefa\ Lemma}  We have: 
\roster
\smallskip\item $R_a$ is
$k$ Osserman for any
$k$ and for any $a$. 
\smallskip\item $R_a$ is Jordan Osserman of type $(r,s)$ if and only
if $\rank R_{\Phi_a}$ is constant on $\Gr_{r,s}(V)$.
\endroster\endproclaim

\demo{Proof}
It is immediate from the definition that $\Phi_a^2=0$. Thus $\range\Phi_a$ is totally
isotropic. We use equation (\brefba) to see that
$\JJ_a(x)y=3(\Phi_ax,y)\Phi_ax$. Consequently 
$$\JJ_a(v_1)\JJ_a(v_2)y=9(\Phi_av_1,y)(\Phi_av_1,\Phi_av_2)\Phi_av_2=0$$ for any vectors $v_1$
and
$v_2$ in $V$. Thus $\JJ_a(\sigma)^2=0$ for any non-degenerate subspace $\sigma$. This implies
that $0$ is the only eigenvalue of $\JJ_a(\sigma)$ and hence $R_a$ is $k$ Osserman for any $k$;
this proves assertion (1). 

Since $\JJ_a(\sigma)^2=0$, the Jordan normal form of $\JJ_a(\sigma)$ is determined by the rank
of $\JJ_a(\sigma)$; assertion (2) now follows directly.
\qed\enddemo

We will use the following Lemma to study the rank of $\JJ_a(\sigma)$.

\setref\drefb
\proclaim{\drefb\ Lemma} Let $v_1$,....,$v_k$ be linearly independent vectors in $V$. Then:
\roster
\smallskip\item There exist vectors $w_1$,...,$w_k$ in $V$ so that $(v_i,w_j)=\delta_{ij}$.
\smallskip\item Let $Ty:=c_1(v_1,y)v_1+...+c_k(v_k,y)v_k$ define a linear transformation of $V$,
where
$c_1$,...,$c_k$ are non-zero constants. Then the rank of $T$ is $k$.
\endroster\endproclaim

\demo{Proof} Let $V^*$ be the associated dual vector space of linear maps
from $V$ to $\real$. We define a linear map $\psi:V\rightarrow V^*$ by setting
$\psi(w)v=(v,w)$. Since the inner product on $V$ is non-degenerate, $\psi$ is injective. Since
$\dim V=\dim V^*$, $\psi$ is bijective. We can extend the collection $\{v_1,...,v_k\}$ to a
basis for $V$; thus without loss of generality, we may assume $k=\dim V$. Let $\{v^1,...,v^k\}$
be the corresponding dual basis for
$V^*$. Set
$w_i:=\psi^{-1}v_i$. Assertion (1) follows as
$$(v_i,w_j)=\psi(w_j)v_i=v^j\cdot v_i=\delta_{ij}.$$

Let $T$ be the transformation of assertion (2). It is clear from the definition that $\range
T\subset\text{span}\{v_1,...,v_k\}$. Since
$(v_i,w_j)=\delta_{ij}$, we have $Tw_j=c_jv_j$. Since $c_j\ne0$, $v_j\in\range T$. It now
follows that $\text{span}\{v_1,...,v_k\}=\range T$. Thus the rank of $T$ is $k$.
\qed\enddemo

We conclude our preparation for the proof of Theorem \brefd\ with the following
Lemma.

\setref\drefc
\proclaim{\drefc\ Lemma} We have: 
\roster
\smallskip\item If $2a<p$ and if $1\le r\le p-1$, then $R_a$ is not Jordan Osserman of type
$(r,s)$ for any $0\le s\le q$.
\item If $2a=p$, if $1\le r\le p-1$, and if $1\le s\le q-1$, then $R_a$ is not Jordan Osserman
of type $(r,s)$.
\smallskip\item $R_a$ is Jordan Osserman of type $(p,0)$.
\smallskip\item If $2a=p$ and if $1\le r\le p-1$, then $R_a$ is Jordan Osserman of type $(r,0)$.
\endroster\endproclaim

\demo{Proof} Let
$\{e_1^-,...,e_p^-,e_1^+,...,e_q^+\}$ be the orthonormal basis for $V$ used in equation
(\CREFbba) to define $\Phi_a$. We have $\JJ_a(e_i^+)=\JJ_a(e_i^-)$; this vanishes if
$i>2a$. Let
$$\tau:=\cases
\pcases{$\Pspan\{e_2^-,...,e_r^-,e_2^+,...,e_s^+\}$}{$r\ge2,\ s\ge2,$}\cr
\pcases{$\Pspan\{e_2^-,...,e_r^-\}$}{$r\ge2,\ s\le1,$}\cr
\pcases{$\Pspan\{e_2^+,...,e_s^+\}$}{$r\le1,\ s\ge2,$}\cr
\pcases{$\{0\}$}{$r\le1,\ s\le1.$}\endcases$$

To prove assertion (1), we will exhibit two subspaces
$\sigma_1$ and $\sigma_2$ of type $(r,s)$ so that
$\rank\{\JJ_a(\sigma_1)\}\ne\rank\{\JJ_a(\sigma_2)\}$. Let
$$\eqalign{
&\sigma_1:=\cases
\pcases{$\tau\oplus\Pspan\{e_1^-,e_1^+\}$\quad\qquad\phantom{..}}{$s\ge1,$}\cr
\pcases{$\tau\oplus\Pspan\{e_1^-\}$}{$s=0$,}\endcases\cr
&\sigma_2:=\cases
\pcases{$\tau\oplus\Pspan\{e_p^-,e_1^+\}$\quad\qquad\phantom{..}}{$s\ge1,$}\cr
\pcases{$\tau\oplus\Pspan\{e_p^-\}$}{$s=0.$}\endcases}$$  
The index `$1$' does not appear among the indices comprising the basis for $\tau$ so there is no
`interaction'. Furthermore, since $2a<p$, $\JJ_a(e_p^-)=0$. We have the cancellation
$(e_1^-,e_1^-)\JJ_a(e_1^-)+(e_1^+,e_1^+)\JJ_a(e_1^+)=0$. We may now use Lemma \drefb\ to
compute:
$$\eqalign{
&\rank\{\Cal{J}_a(\sigma_1)\}=\cases
\pcases{$\rank\{\Cal{J}_a(\tau)\}$}{$s\ge1,$}\cr
\pcases{$\rank\{\Cal{J}_a(\tau)\}+1$}{$s=0,$}\endcases\cr
&\rank\{\Cal{J}_a(\sigma_2)\}=\cases
\pcases{$\rank\{\Cal{J}_a(\tau)\}+1$}{$s\ge1,$}\cr
\pcases{$\rank\{\Cal{J}_a(\tau)\}$}{$s=0.$}\endcases}$$
This shows that $\JJ_a(\sigma_1)$ and $\JJ_a(\sigma_2)$ have different ranks and thereby
completes the proof of assertion (1).

Let $2a=p$, let $1\le r\le p-1$, and let $1\le s\le q-1$. As in the proof of assertion (1), we
will construct two subspaces $\sigma_1$ and $\sigma_2$ of type $(r,s)$ so that
$\JJ_a(\sigma_1)$ and $\JJ_a(\sigma_2)$ have different ranks. We define 
$$\eqalign{
&\sigma_1:=\phantom{\quad}\tau\oplus\Pspan\{e_1^-,e_1^+\},\phantom{a..}\text{ and}\cr
&\sigma_2:=\cases
\pcases{$\tau\oplus\Pspan\{e_{r+1}^-,e_1^+\}$}{$r\ge s,$}\cr
\pcases{$\tau\oplus\Pspan\{e_1^-,e_{s+1}^+\}$}{$r<s.$}\endcases}$$
Again, note that the index `1' does not appear among the indices comprising the basis for
$\tau$. If
$r\ge s$ (resp. $r<s$), then the index $r+1$ (resp. $s+1$) does not appear among these indices
either. If
$s+1\le p$, then $\JJ_a(e_{s+1}^+)\ne0$; if $s+1>p$, then $\JJ_a(e_{s+1}^-)=0$. Since
$r+1\le p=2a$, $\JJ_a(e_{r+1}^-)\ne0$. Thus:
$$\eqalign{
&\rank\{\Cal{J}_a(\sigma_1)\}=\phantom{\quad}\rank\{\Cal{J}_a(\tau)\}\cr
&\rank\{\Cal{J}_a(\sigma_2)\}=\cases
\pcases{$\rank\{\Cal{J}_a(\tau)\}+2$}{$r\ge s$}\cr
\pcases{$\rank\{\Cal{J}_a(\tau)\}+2$}{$r<s<p,$}\cr
\pcases{$\rank\{\Cal{J}_a(\tau)\}+1$}{$r<p\le s.$}\endcases}$$
Since $\JJ_a(\sigma_1)$ and $\JJ_a(\sigma_2)$ have different ranks, $R_a$ is not Jordan
Osserman of type $(r,s)$. This establishes assertion (2).

If we can show that $\rank\{\JJ_a(\sigma)\}=2a$ for every maximal timelike subspace
$\sigma$ of
$V$, we may then use Lemma \drefa\ to show that $R_a$ is Jordan Osserman of type $(p,0)$
which will prove assertion (3). Since
$$\range\{\JJ_a(\sigma)\}\subset\{\Phi_ae_1^-,...,\Phi_ae_{2a}^-\},\qquad
\rank\{\JJ_a(\sigma)\}\le2a.$$
We suppose that $\rank\{\JJ_a(\sigma)\}<2a$ and argue for a contradiction. Let
$$W:=\text{span}\{e_1^-,...,e_{2a}^-\}.$$
As $\dim W=2a$ and as $\rank\{\JJ_a(\sigma)\}<2a$, we have
$\ker\{\JJ_a(\sigma)\}\cap W\ne\{0\}$. Thus we may choose $0\ne w\in W$ with
$\JJ_a(\sigma)w=0$. Let $\{v_1,...,v_p\}$ be an orthonormal basis for $\sigma$. We use
equations (\arefac) and (\brefba) to compute:
$$0=(\JJ_a(\sigma)w,w)=-3(\Phi_av_1,w)(\Phi_av_1,w)-...-3(\Phi_av_p,w)(\Phi_av_p,w).$$
This implies that 
$$0=(\Phi_av_i,w)=-(v_i,\Phi_aw)\quad\text{ for }1\le i\le p.$$
Consequently,
$\Phi_aw\perp\sigma$. Since $\sigma$ is a maximal timelike subspace, $\Phi_aw$ either
vanishes or is spacelike. Since $\range\Phi_a$ is totally isotropic, $\Phi_aw$ is a null vector,
not a spacelike vector. Thus we must have that
$\Phi_aw=0$. This is false as
$\Phi_a$ is injective on $W$. This contradiction shows that $\rank\{\JJ_a(\sigma)\}=2a$ and
hence
$R_a$ is Jordan Osserman of type
$(p,0)$;
assertion (3) is established.

Finally, suppose that $2a=p$. Let $\sigma\in\Gr_{r,0}(V)$. If we can show that the rank of
$\JJ_a(\sigma)$ is $r$, then it would follow by Lemma \drefa\ that $R_a$ is Jordan
Osserman of type $(r,0)$. Let $\{v_1,...,v_r\}$ be an orthonormal basis for $\sigma$. We use
equations (\arefac) and (\brefba) to see that
$$\JJ_a(\sigma)y=-3\{(y,\Phi_a(v_1))\Phi_a(v_1)+...+(y,\Phi_a(v_r))\Phi_a(v_r)\}.$$
Thus by Lemma \drefb, if $\{\Phi_a(v_1),...,\Phi_a(v_r)\}$ is a linearly independent
set, then $\rank\JJ_a(\sigma)=r$. We suppose the contrary, i.e. that the set
$\{\Phi_a(v_1),...,\Phi_a(v_r)\}$ is linearly dependent, and argue for a
contradiction. Choose $0\ne v\in\sigma$ so that $\Phi_a(v)=0$. We expand:
$$\eqalign{
v&=c_1^-e_1^-+...+c_p^-e_1^-+c_1^+e_1^++...+c_q^+e_q^+\cr
\Phi_av&=(c_1^+-c_1^-)\Phi_a(e_1^+)+...+(c_p^+-c_p^-)\Phi_a(e_p^+).
}$$
Since $\Phi_av=0$, we have $c_i^+=c_i^-$ for $1\le i\le p$. Thus
$$\eqalign{
(v,v)&=(c_1^+)^2-(c_1^-)^2+...+(c_p^+)^2-(c_p^-)^2+(c_{p+1}^+)^2+...+(c_q^+)^2\cr
     &=(c_{p+1}^+)^2+...+(c_q^+)^2\ge0.}$$
Since $v\ne0$ and since $\sigma$ is timelike, $(v,v)<0$. This contradiction shows that
$\{\Phi_av_1,...,\Phi_av_r\}$ is a linearly independent set. Thus
$\rank\{\JJ_a(\sigma)\}=r$ and hence $R_a$ is Jordan Osserman of type $(r,0)$. This completes
the proof of the last assertion of the Lemma.
\qed\enddemo

\setref\drefe
\definition{\drefe\ Remark} We note that we can interchange the roles of spacelike and
timelike vectors in Lemma \drefc\ by changing the sign of the inner product.
\enddefinition

\demo{Proof of Theorem {\rm \brefd}} 
Suppose that $2a<p\le q$. We use Lemma \drefc\ (3) to see that $R_a$ is Jordan Osserman
of type $(p,0)$. Then dually by Theorem \brefb, we have that $R_a$ is Jordan Osserman
of type $(0,q)$. By Lemma \drefc\ (1), since $2a<p$, $R_a$ is not Jordan Osserman of type
$(r,s)$ for $1\le r\le p-1$. Similarly, by Lemma \drefc\ (1) and Remark \drefe, since $2a<q$,
$R_a$ is not Jordan Osserman of type $(r,s)$ for $1\le s\le q-1$. Assertion (1) now follows.

Suppose that $2a=p<q$. We use Lemma \drefc\ (4) to see that $R_a$ is Jordan Osserman of
type $(r,0)$ if $1\le r\le p-1$. Dually, by Theorem \brefb, $R_a$ is also Jordan
Osserman of type $(p-r,q)$. By Lemma \drefc\ (1) and Remark \drefe, since $2a<q$, $R_a$ is not
Jordan Osserman of type $(r,s)$ if $1\le s\le q-1$. This establishes assertion (2).

Suppose finally that $2a=p=q$. Since $2a=p$, we use Lemma \drefc\ (4) to see that $R_a$
is Jordan Osserman of type $(r,0)$ if $1\le r\le p-1$. Since $2a=q$, similarly we have, by
Remark \drefe, that
$R_a$ is Jordan Osserman of type $(0,s)$ if $1\le s\le q-1$. The remaining values $(p-r,q)$ and
$(p,q-s)$ then follow dually by Theorem \brefb. If $1\le r\le p-1$ and $1\le s\le
q-1$, then $R_a$ is not Jordan Osserman by Lemma \drefc\ (2). All the assertions of
Theorem \brefd\ are now proved.
\qed\enddemo

\enddocument
\bye